\documentclass{amsart}

\usepackage{etex}

\usepackage{hyperref}
\usepackage{txfonts, amsmath,amstext,amsthm,amscd,amsopn,verbatim,amssymb, amsfonts, amsrefs}
\usepackage{fullpage}

\usepackage[bbgreekl]{mathbbol}

\usepackage{tikz}

\usetikzlibrary{matrix}
\usetikzlibrary{shapes}
\usetikzlibrary{arrows}
\usetikzlibrary{calc,3d}
\usetikzlibrary{decorations,decorations.pathmorphing}
\usetikzlibrary{through}
\tikzset{ext/.style={circle, draw,inner sep=1pt},int/.style={circle,draw,fill,inner sep=1pt},nil/.style={inner sep=1pt}}
\tikzset{exte/.style={circle, draw,inner sep=3pt},inte/.style={circle,draw,fill,inner sep=3pt}}
\tikzset{diagram/.style={matrix of math nodes, row sep=3em, column sep=2.5em, text height=1.5ex, text depth=0.25ex}}
\tikzset{diagram2/.style={matrix of math nodes, row sep=0.5em, column sep=0.5em, text height=1.5ex, text depth=0.25ex}}

\usepackage{tikz-cd}

\theoremstyle{plain}

\newtheorem{thm}{Theorem}[section]

\newtheorem{prop}[thm]{Proposition}
\newtheorem{cor}[thm]{Corollary}
\newtheorem{lemm}[thm]{Lemma}

\theoremstyle{definition}

\newcommand{\alg}[1]{\mathfrak{{#1}}}

\newcommand{\p}{\partial}

\newcommand{\R}{{\mathbb{R}}}

\newcommand{\K}{{\mathbb{K}}}
\newcommand{\Q}{{\mathbb{Q}}}

\newcommand{\HGC}{{\mathrm{HGC}}}

\newcommand{\Poiss}{\mathsf{Pois}}

\newcommand{\Com}{\mathsf{Com}}

\newcommand{\bpm}{\begin{pmatrix}}
\newcommand{\epm}{\end{pmatrix}}

\newcommand{\GC}{\mathrm{GC}}

\newcommand{\tmid}{|}

\newcommand{\hotimes}{\mathbin{\hat\otimes}}

\newcommand{\e}{\mathsf{e}}

\DeclareMathOperator{\Emb}{Emb}
\DeclareMathOperator{\Imm}{Imm}
\DeclareMathOperator{\Embbar}{\overline{Emb}}

\newcommand{\lD}{\mathsf{D}}

\newcommand{\Map}{\mathrm{Map}}

\newcommand{\Aut}{\mathrm{Aut}}

\newcommand{\bS}{\mathsf{N}}

\newcommand{\stp}{{{}^*\mathfrak{p}}}
\newcommand{\pp}{{\hat{\alg p}}}
\newcommand{\Npp}{{\bS\hat{\alg p}}}

\newcommand{\beq}[1]{\begin{equation}\label{#1}}
\newcommand{\eeq}{\end{equation}}

\DeclareMathOperator*{\holim}{holim}

\newcommand{\ar}{\mathrm{ar}}

\begin{document}
\title{On truncated mapping spaces of $E_n$ operads}

%

\author{Thomas Willwacher}
\address{Department of Mathematics \\ ETH Zurich \\
R\"amistrasse 101 \\
8092 Zurich, Switzerland}
\email{thomas.willwacher@math.ethz.ch}

\date{January 30, 2017}

\thanks{
T.W. has been partially supported by the NCCR SwissMAP funded by the Swiss National Science Foundation, and the ERC starting grant 678156 GRAPHCPX}


\begin{abstract}
The mapping spaces of the $r$-truncated versions of the $E_n$ operads appear as the $r$-th stage of the Taylor tower for long embedding spaces.
It has been shown in \cite{FTW} that their rational homotopy groups can be expressed through graph homology in the limit $r\to \infty$.
For finite $r$ only a part of the graph homology appears in the homotopy groups of the Taylor tower at this stage, possibly along with some additional unstable homotopy groups.
In this paper we study the convergence properties in some more detail. In particular, we provide bounds for the stage $r$ at which the various graph homology classes start appearing, and we provide degree bounds for the unstable (i.e., vanishing as $r\to \infty$) terms in the homotopy groups of the tower.
\end{abstract}

\maketitle


\section{Introduction }
The spaces $\Emb_\p(\R^m,\R^n)$ of long embeddings of $\R^m$ in $\R^n$, i.e., embeddings fixed outside of a compact to be the standard embedding, are well studied objects of algebraic topology.
In particular, it is known that they can be accessed using the Goodwillie-Weiss embedding calculus \cite{GW}.
More precisely, we consider the homotopy fiber $\Embbar_\p(\R^m,\R^n)$ of the inclusion $\Emb_\p(\R^m,\R^n)\to \Imm_\p(\R^m,\R^n)$ of embeddings into immersions.
Then the Taylor tower converges for $n-m\geq 3$,
\[
\Embbar_\p(\R^m,\R^n) \simeq 
T_\infty\Embbar_\p(\R^m,\R^n)
:= \holim_k T_k\Embbar_\p(\R^m,\R^n).
\]
Furthermore, the $k$-th stage of the Taylor tower $T_k\Embbar_\p(\R^m,\R^n)$ can be expressed as the $m+1$-st loop space of the derived mapping spaces of the $k$-truncated little disks operads, and one has an equivalence of towers \cite{WBdB2,Weiss2,DucT}
\begin{equation}\label{eq:towers}
\begin{tikzcd}
*\ar{d}{\simeq} &
\Omega^{m+1}\Map_{\leq 1}^h(\lD_m,\lD_n)\ar{d}{\simeq}\ar{l} &
\cdots\ar{l} &
\Omega^{m+1}\Map_{\leq k}^h(\lD_m,\lD_n)\ar{d}{\simeq}\ar{l} &
\cdots\ar{l} &
\Omega^{m+1}\Map^h(\lD_m,\lD_n)\ar{d}{\simeq}\ar{l}\\
T_0\Embbar_{\partial}(\R^m,\R^n) &
T_1\Embbar_{\partial}(\R^m,\R^n)\ar{l} &
\cdots\ar{l} &
T_k\Embbar_{\partial}(\R^m,\R^n)\ar{l} &
\cdots\ar{l} &
T_\infty\Embbar_{\partial}(\R^m,\R^n)\ar{l}
\end{tikzcd},
\end{equation}

Here a $k$-truncated operad is similar data as an operad, except that one forgets about all operations of arities $>k$, and all compositions involving such operations.

The rational homotopy type of the mapping space $\Map^h(\lD_m,\lD_n)$ has been determined in \cite{FTW}. It has been shown there that an $L_\infty$-algebra model of this space is provided by the hairy graph complex $\HGC_{m,n}$.
We refer the reader to \cite[section 7]{FTW} or \cite{Turchin3} for the definition and recollections about the hairy graph complex, and just remark here that its elements are formal $\Q$-linear series of graphs with external legs ("hairs") such as the following.
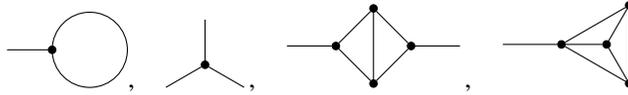
\[
\begin{tikzpicture}[scale=.5]
\draw (0,0) circle (1);
\draw (-180:1) node[int]{} -- +(-1.2,0);
\end{tikzpicture},
\quad\begin{tikzpicture}[scale=.6]
\node[int] (v) at (0,0){};
\draw (v) -- +(90:1) (v) -- ++(210:1) (v) -- ++(-30:1);
\end{tikzpicture}
{\,},
\quad\begin{tikzpicture}[scale=.5]
\node[int] (v1) at (-1,0){};\node[int] (v2) at (0,1){};\node[int] (v3) at (1,0){};\node[int] (v4) at (0,-1){};
\draw (v1)  edge (v2) edge (v4) -- +(-1.3,0) (v2) edge (v4) (v3) edge (v2) edge (v4) -- +(1.3,0);
\end{tikzpicture}
{\,},
\quad\begin{tikzpicture}[scale=.6]
\node[int] (v1) at (0,0){};\node[int] (v2) at (180:1){};\node[int] (v3) at (60:1){};\node[int] (v4) at (-60:1){};
\draw (v1) edge (v2) edge (v3) edge (v4) (v2)edge (v3) edge (v4)  -- +(180:1.3) (v3)edge (v4);
\end{tikzpicture}
\]
In particular, the hairy graph homology appears as the rational homotopy groups of the embedding spaces $\Embbar_\p(\R^m,\R^n)$ for $n-m\geq 3$, as has been shown in \cite{FTW}, and earlier in \cite{Turchin2} for a weaker degree range.
In other words, the stable part of the rational homotopy groups of the Taylor tower is given by the hairy graph homology.
The purpose of the present paper is to study in more detail the precise convergence properties of the Taylor tower.
To roughly summarize our results, we have the following stabilization behavior.
The part of the hairy graph homology with $g$ loops (i.e., the first Betti number of graphs involved is $g$) and $h$ hairs, $H(\HGC_{m,n}^{g\text{-loop},h\text{-hairs}})$ starts appearing in the Taylor tower roughly at stage $k=g+h$, and survives as $k\to \infty$. Apart from this stable piece there are (generally) "spurious" classes in a certain degree range that can be bounded.
These classes are short-lived, in that they vanish in the immediate next stage after the stage they first appeared in.

Let us remark that the objects we study will be the derived mapping spaces $\Map^h_{\leq k}(\lD_m,\lD_n^\Q)$ between the little $m$-disks operad $\lD_m$, and the rationalization $\lD_n^\Q$. It has been shown in \cite{FTW} that the natural map of towers
\[
\begin{tikzcd}
\cdots &
\Map^h_{\leq k}(\lD_m,\lD_n) \ar{l} \ar{d}{\simeq_{\Q}}& 
\Map^h_{\leq k+1}(\lD_m,\lD_n) \ar{l} \ar{d}{\simeq_{\Q}}&
\cdots  \ar{l} \\
\cdots &
\Map^h_{\leq k}(\lD_m,\lD_n^\Q) \ar{l} &
\Map^h_{\leq k+1}(\lD_m,\lD_n^\Q) \ar{l} & 
\cdots \ar{l}
\end{tikzcd}
\]
is a (termwise) rational homotopy equivalence for $n-m\geq 3$. However, even for $n=m$ one can study the derived mapping spaces
\[
\Map^h_{\leq k}(\lD_n,\lD_n^\Q).
\]
In the stable case $k\to\infty$ their rational homotopy groups relative to the basepoint given by the rationalization map $\lD_n\to \lD_n^\Q$ are computed in \cite{FTW}, and can be expressed through the homology of the Kontsevich graph complexes $\GC_n^2$. These latter grpah complexes are defined similarly to $\HGC_{m,n}$ above, except that the graph do not have hairs.
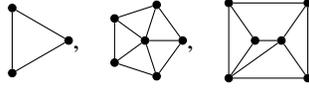
\[
\begin{tikzpicture}[baseline=-.65ex, scale=.5]
\node[int] (v1) at (0:1) {};
\node[int] (v2) at (120:1) {};
\node[int] (v3) at (240:1) {};
\draw (v1) edge (v2) edge (v3) (v3) edge (v2) ;
\end{tikzpicture},
\quad\begin{tikzpicture}[baseline=-.65ex, scale=.5]
\node[int] (c) at (0,0){};
\node[int] (v1) at (0:1) {};
\node[int] (v2) at (72:1) {};
\node[int] (v3) at (144:1) {};
\node[int] (v4) at (216:1) {};
\node[int] (v5) at (-72:1) {};
\draw (v1) edge (v2) edge (v5) (v3) edge (v2) edge (v4) (v4) edge (v5)
(c) edge (v1) edge (v2) edge (v3) edge (v4) (c) edge (v5);
\end{tikzpicture},
\quad\begin{tikzpicture}[baseline=-.65ex, scale=.5]
\node[int] (c) at (0.7,0){};
\node[int] (v1) at (0,-1) {};
\node[int] (v2) at (0,1) {};
\node[int] (v3) at (2.1,-1) {};
\node[int] (v4) at (2.1,1) {};
\node[int] (d) at (1.4,0) {};
\draw (v1) edge (v2) edge (v3)  edge (d) edge (c) (v2) edge (v4) edge (c) (v4) edge (d) edge (v3) (v3) edge (d) (c) edge (d);
\end{tikzpicture}
\]
In this situation we derive similar (albeit a bit weaker) results for the stabilization behavior. The $g$-loop graph homology $H(\GC^{2,g\text{-loop}}_{n})$ starts appearing in the tower roughly at stage $k=g$. Apart from this stable piece, there are generally unstable classes in a certain range.


\subsection{Statement of results}

As a corollary of the discussion in section 10 (see in particular Theorem 10.12 and Corollary 10.13) of \cite{FTW} we can in particular show the following result.
For $n\geq m\geq 1$ and $n\geq 2$ (and $*\geq 1$) consider the map defined in \cite{FTW}
\begin{equation}\label{equ:mapmn}
I_k: H_{*-1}(\HGC_{m,n})\cong \pi_*(\Map^h(\lD_m, \lD_n^{\Q}), *) \to \pi_*(\Map^h_{\leq k}(\lD_m, \lD_n^{\Q}), *)
\end{equation}
where the basepoint $*$ is taken to be the trivial map $\lD_m\to | \Com^c| \to \lD_n^{\Q}$.

\begin{thm}\label{thm:mainHGC}
The kernel of the map $I_k$ of \eqref{equ:mapmn} is, in degrees $*\geq 1$,
\[
\bigoplus_{g+h>k} H_{<(n-m-2)(g+h-1)+mg+h-k+1}(\HGC_{m,n}^{g\text{-loop},h\text{-hairs}}) 
\subset H_*(\HGC_{m,n}).
\]
The cokernel of the map $I_k$ can be expressed as 
\begin{equation}\label{equ:HGCcokernel}
\bigoplus_{g\geq 1,h\geq 1 \atop h\leq k \leq g+h-1}
V_{g,h,k},
\end{equation}
where $V_{g,h,k}$ is a finite dimensional graded vector space concentrated in the single degree 
\[
(n-m-2)(g+h-1)+mg+h-k+2. 
\]
\end{thm}
The statement of the Theorem is that the arity truncation of the mapping spaces corresponds approximately to the truncation by the number $g+h$ (loops plus hairs) on the hairy graph complex side. 
We do not exclude that $V_{g,h,k}$ is zero in some or many cases.
However, it will be clear from the proof that the classes in $V_{g,h,k}$ already vanish at the $k+1$-st stage of the tower.

Also note that one has the following degree bounds on $H(\HGC_{m,n}^{g\text{-loop},h\text{-hairs}})$.
\begin{lemm}\label{lem:degboundHGC}
The hairy graph homology
$
H_d(\HGC_{m,n}^{g\text{-loop},h\text{-hairs}})
$
is concentrated in the range (i.e., is zero outside the range)
\[
(n-m-2)(g+h-1)+(m-1)g +1 \leq d \leq (n-m-2)(g+h-1)+mg +1.
\]
\end{lemm}

Similarly, for $n\geq 2$ (and $*\geq 1$) we consider the map defined in \cite{FTW}
\begin{equation}\label{equ:mapnn}
J: H_{*}(\GC^2_{n})\cong \pi_*(\Map^h(\lD_n, \lD_n^{\Q}),*_{id}) \to \pi_*(\Map^h_{\leq k}(\lD_n, \lD_n^{\Q}), *_{id}),
\end{equation}
where as the basepoint we use the canonical map $\lD_n\to \lD_n^\Q$.
We also note that the corresponding component of the derived mapping space is weakly equivalent to the identity component of the homotopy automorphism space $\Aut^h(\lD_n^{\Q})$, so that we may equivalently consider the maps 
\[
J': H_{*}(\GC^2_{n})\cong \pi_*(\Aut^h(\lD_n^{\Q}),*_{id}) \to \pi_*(\Aut^h_{\leq k}(\lD_n^{\Q}), *_{id}).
\]

In this situation one has the following result.
\begin{prop}\label{prop:mainGC}
The restriction of $J$ of \eqref{equ:mapnn} to the $g$-loop part
\[
H_{d}(\GC^{2,g\text{-loop}}_{n})\to \pi_d(\Map^h_{\leq k}(\lD_n, \lD_n^{\Q}), *_{id})
\]
is an injection in positive degrees $d$ satisfying
\[
d\geq (n-2) g +4- 2k +\min(k-1,g) =: d_{g,k,1}
\]
and is zero for positive $d$ such that
\[
d\leq (n-2)g +3- 2k=:d_{g,k,0}.
\]
Furthermore, the cokernel of the map $J$ has the form
\[
\bigoplus_{g} W_{g,k},
\]
where $W_{g,k}$ is a finite dimensional graded vector space concentrated in degrees $d$ with
\[
d_{g,k,0} < d \leq d_{g,k,1}.
\]
\end{prop}
We also remark on the following well known degree bounds for the graph complex.
\begin{lemm}[\cite{Will}]\label{lem:degboundGC}
For $g\geq 2$ the graph homology $H_{d}(\GC^{2,g\text{-loop}}_{n})$ vanishes outside of the range
\[
(n-3)g +3 \leq d \leq (n-2)g+1.
\]
For $g\geq 3$ the upper bound can be improved by one to $(n-2)g$.
\end{lemm}
We also note that the loop order $\leq 2$-pieces of the graph cohomology $H(\GC_n^2)$ are spanned by the classes represented by the following graphs.
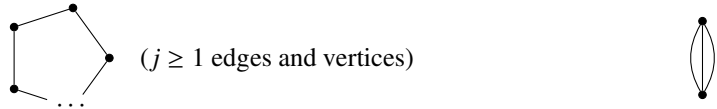
\begin{align*}
&
\begin{tikzpicture}[baseline=-.65ex, scale=.7]
\node[int] (v1) at (0:1) {};
\node[int] (v2) at (72:1) {};
\node[int] (v3) at (144:1) {};
\node[int] (v4) at (216:1) {};
\node (v5) at (-72:1) {$\cdots$};
\draw (v1) edge (v2) edge (v5) (v3) edge (v2) edge (v4) (v4) edge (v5);
\end{tikzpicture}
\quad\text{($j\geq 1$ edges and vertices) }&
&
\begin{tikzpicture}[baseline=-.65ex, scale=.7]
\node[int] (v1) at (0,-.7) {};
\node[int] (v2) at (0,.7) {};
\draw (v1) edge (v2) edge[bend left] (v2)  edge[bend right] (v2);
\end{tikzpicture}
\end{align*}
(Mind that some of them are 0, depending on the parity of $n$.)

Combining Lemma \ref{lem:degboundGC} with Proposition \ref{prop:mainGC} we find in particular.
\begin{cor}
	The restriction of the map $J$ to $H(\GC^{2,g\text{-loop}}_{n})$ is injective in positive degrees for $k\geq g\geq 2$.
\end{cor}

Roughly one can hence say that the truncation by arity of the automorphisms of the little disks operad corresponds to truncation by loop order in the graph complex. Mind that Proposition \ref{prop:mainGC} however leaves a "gap" of degrees between $d_{g,k,0}$ and $d_{g,k,1}$, within which we do not understand the behaviour of the map $J$ well.
Furthermore, the situation is a bit different in loop order one. One can check that the loop order $1$ graph with $j$ edges appears on the $\frac{j+3}{2}$-th stage of the tower, by hand computation.

\section*{Acknowledgements}
This paper is a short addendum to the paper \cite{FTW} by B. Fresse, V. Turchin and the author.
In fact, all results presented here are mere degree counting exercises given that paper. We compile them here for easier reference, and since they are not so easy to extract from \cite{FTW}.
The author is very grateful to B. Fresse and V. Turchin for many fruitful discussions. They could also be mentioned as coauthors, but preferred not to assume this role.


\section{Proofs for the "hairy" case}
It is shown in \cite{FTW} that the (positive) homotopy groups of the mapping spaces $\Map^h(\lD_n, \lD_n^{\Q})$ are computed as the homology of the Hopf cooperadic deformation complexes
\[
K(\stp_n, \Poiss_m\{m\}).
\]
For the definition of these complexes we refer the reader to \cite[section 3.4]{FTW}.
For us only the following facts will be important:
\begin{itemize}
	\item As a graded vector space we have 
	\begin{equation}\label{equ:Kdef}
	K(\stp_n, \Poiss_m\{m\})
	= \prod_{r=2}^\infty \Npp_n(r)[-1]\hotimes_{S_r} \Poiss_m\{m\}(r),
	\end{equation}
	where $\pp_n(r)$ is the (higher dimensional) Drinfeld-Kohno Lie algebra generated by symbols $t_{ij}$ of degree $n-2$ with $1\leq i,j\leq r$, and relations as in \cite[Theorem 4.6]{FTW}. The subspace $\Npp_n(r)\subset \pp_n(r)$ is spanned by those Lie monomials in which each index $i=1,\dots, r$ occurs at least once.
	It will be important to us that $\Npp_n(r)$ is naturally graded, each generator carrying degree $+1$.
	The operad $\Poiss_m$ (the $m$-Poisson operad) is generated by two binary operations, the commutative product of degree $0$, and the Lie bracket of degree $m-1$. It carries a natural grading by the number of brackets.
	Mind that in \cite{FTW} $\Poiss_m$ is replaced by $\e_m:=H(\lD_m)$, which is $\Poiss_m$ for $m\geq 2$ and the associative operad for $m=1$. However, even for $m=1$ the complexes produced are quasi-isomorphic, so it makes no difference for our paper. 
	\item The differential raises the arity $r$ by one, but preserves the grading by number of generators on the $\pp_n'(r)$.
	It furthermore raises the number of brackets in $\Poiss_m\{m\}$ by one.
	This means in particular that the following two numbers are invariant under the differential
	\begin{align*}
	h &:= r-(\text{brackets})\\
	g &:= (\text{generators})-h+1 .
	\end{align*}
	We call $g$ the loop order and $h$ the number of hairs.
	The number of hairs obviously satisfies $h\geq 1$.
	Furthermore, the loop order is always non-negative, $g\geq 0$, since the number of generators must be at least $r-1$.
	Our deformation complex splits into a direct product of subcomplexes $C_{m,n,g,h}\subset K(\stp_n, \Poiss_m\{m\})$ according to these numbers,
	\[
	K(\stp_n, \Poiss_m\{m\})
	= \prod_{g,h} C_{m,n,g,h}.
	\]
	
	\item As shown in \cite{FTW} one has, for all $d$, an identification 
	\[
	H_d(K(\stp_n, \Poiss_m\{m\})) \cong H_d(\HGC_{m,n})
	\]
	in such a way that the grading by $g$ and $h$ above is precisely mapped to grading on the hairy graph homology by the loop number $g$ and hair number $h$ of graphs.
	Furthermore, for positive $d$ one has  \cite{FTW}
	\begin{equation}\label{equ:Kpi}
	H_{d-1}(K(\stp_n, \Poiss_m\{m\})) \cong 
	\pi_d(\Map^h(\lD_m, \lD_n^{\Q}), *) .
	\end{equation}
\end{itemize}

Similarly, one can consider a truncated version of the above story.
In particular, it is shown in \cite[section 10.3]{FTW} that the homotopy groups of the mapping spaces of truncated operads $\Map^h_{\leq k}(\lD_m, \lD_n^{\Q})$ are computed as the homology of the Hopf cooperadic deformation complexes
\[
K(\stp_n\tmid_{\leq k}, \Poiss_m\{m\}\tmid_{\leq k})
\]
which can defined as the quotients (cf. \eqref{equ:Kdef})
\[
K(\stp_n\tmid_{\leq k}, \Poiss_m\{m\}\tmid_{\leq k})
=
K(\stp_n, \Poiss_m\{m\})
/ \prod_{r\geq k+1} \Npp_n(r)\hotimes \Poiss_m\{m\}(r).
\]
Concretely, one then has in positive degree $d$
\begin{equation}\label{equ:Kpitrunc}
H_{d-1}(K(\stp_n\tmid_{\leq k}, \Poiss_m\{m\}\tmid_{\leq k})) \cong 
\pi_d(\Map^h_{\leq k}(\lD_m, \lD_n^{\Q}), *) .
\end{equation}
Similarly to $C_{m,n,g,h}$ above we define the quotients 
\[
C_{m,n,g,h,k} = C_{m,n,g,h} / \left(C_{m,n,g,h}\cap  \prod_{r\geq k+1} \Npp_n(r)\hotimes \Poiss_m\{m\}(r)\right)
\]
so that 
\[
K(\stp_n\tmid_{\leq k}, \Poiss_m\{m\}\tmid_{\leq k})\cong \prod_{g,h}C_{m,n,g,h,k}.
\]

\subsection{Proof of Lemma \ref{lem:degboundHGC}}
The cohomological degree $d$ in the deformation complex \eqref{equ:Kdef} can be expressed by the numbers $g,h,r$ as follows:
\begin{equation}\label{equ:d}
\begin{aligned}
d &= 
(n-2)(\text{generators})
+
(m-1)(\text{brackets})
-m(r-1)
+1
\\&=
(n-m-2)(\text{generators})
+
m(\text{generators}+\text{brackets}-r+1)
-\text{brackets}
+1
\\&=
(n-m-2)(g+h-1)
+
mg
+h-r
+1
\end{aligned}
\end{equation}
For fixed $g,h$ the number $r$ can be bounded as follows
\begin{equation}\label{equ:rbounds}
\begin{aligned}
r&\geq \max(h,2)\geq h\\
r&\leq (\text{generators})+1 = h+g.
\end{aligned}
\end{equation}
But then Lemma \ref{lem:degboundHGC} immediately follows.
\hfill\qed

\subsection{Proof of Theorem \ref{thm:mainHGC}}\label{sec:HGCproof}
We fix $g,h$ and we want to study the homology map induced by the quotient map, which models the $g,h$-part of \eqref{equ:mapmn} (see \eqref{equ:Kpi}, \eqref{equ:Kpitrunc})
\[
C_{m,n,g,h}\to C_{m,n,g,h,k}
\]
from the subcomplex of the deformation complex of fixed $g$ and $h$ to its truncated version.
Given the formula for the homological degree of the previous subsection we see immediately that the arity truncation at arity $k$ is the same as the degree truncation at degree 
\[
d_c=(n-m-2)(g+h-1)
+
mg
+h-k
+1.
\]
Concretely, the quotient $C_{m,n,g,h,k}$ is obtained from $C_{m,n,g,h}$ by setting to zero all elements of homological degrees $\leq d_c-1$.
It immediately follows that 
\[
H_d(C_{m,n,g,h,k})
=
\begin{cases}
H_d(C_{m,n,g,h})=H_d(\HGC_{m,n}^{g\text{-loop},h\text{-hair}})& \text{if $d>d_c$} \\
0 &\text{if $d<d_c$} \\
H_d(C_{m,n,g,h})\oplus V &\text{if $d=d_c$}
\end{cases},
\]
where $V$ is some (non-canonical) quotient of the degree $d_c$ piece of $C_{m,n,g,h}$.
We define $V_{g,h,k}$ to be the piece corresponding to $V$ of the cokernel of \eqref{equ:mapmn}, concentrated in degree $d_c+1$.
(The degree shift is due to the shift in \eqref{equ:Kpi}.)

From these observations, Theorem \ref{thm:mainHGC} essentially follows immediately, except for the following cosmetic observation.
First note that the complex $V$ above is zero if $d_c$ is either the bottom degree in the complex $C_{m,n,g,h}$, or a degree in which the complex is zero.
But from the discussion of the previous subsection we know that $C_{m,n,g,h}$ is concentrated in homological degrees $d$ with 
\[
(n-m-2)(g+h-1)
+
(m-1)g +1
\leq d\leq
(n-m-2)(g+h-1)
+ mg +1.
\]
In particular, we can conclude that $V=0$ unless 
\[
(n-m-2)(g+h-1)
+
(m-1)g +1
< d_c\leq
(n-m-2)(g+h-1)
+ mg +1.
\]
Inserting the expression for $d_c$ we obtain the following necessary conditions for $V$ to be nontrivial
\[
-g < h-k \leq 0.
\]
This is reflected in the index ranges in the direct sum \eqref{equ:HGCcokernel} in the Theorem.
\hfill\qed

\section{Proofs for non-hairy case}
We next consider the case $m=n$, and restrict our attention to the component of the canonical map $*_{id}:\lD_n\to \lD_n^{\Q}$ in the mapping space $\Map^h(\lD_n, \lD_n^{\Q})$.
According to \cite{FTW} we can compute the homotopy groups of this connected component by twisting the (relevant) deformation complex with an appropriate Maurer-Cartan element.
Eventually, it turns out that we can compute the homotopy groups as the homology of the complex 
\begin{equation}\label{equ:Ktwisted}
(K(\stp_n, \Poiss_n\{n\}), d_0+d_1),
\end{equation}
already considered above, but with an additional piece $d_1$ contributing to the differential.
For our purposes we only need to know the following two facts from \cite{FTW}:
\begin{itemize}
	\item The piece $d_1$ leaves the arity $r$ invariant, but adds a generator and removes a bracket. This means that the grading by loop number $g$ is still a grading on the complex \eqref{equ:Ktwisted}, but the number of hairs $h$ no longer is.
	Accordingly, our complex splits into a direct product of subcomplexes
	\[
	(K(\stp_n, \Poiss_n\{n\}), d_0+d_1)
	\cong 
	\prod_g D_{n,g}.
	\]
	\item We have that
	\begin{equation}\label{equ:DngGC}
	H_d(D_{n,g}) = H_{d+1}(\GC_n^{2,g\text{-loop}})
	\end{equation}
	for $g\geq 1$ and all $d$, while $H_d(D_{n,0}) = H_{d+1}(\K\ltimes \GC_n^{2,0\text{-loop}})$. Furthermore, for positive $d$
	\[
	H_{d-1}(K(\stp_n, \Poiss_n\{n\}), d_0+d_1)
	\cong
	\pi_d(\Map^h_{\leq k}(\lD_n, \lD_n^{\Q}), *_{id}).
	\]
	
\end{itemize}

Similarly, the homotopy groups for the truncated mapping space $\Map^h_{\leq k}(\lD_n, \lD_n^{\Q})$ are computed by the arity-truncated version of the deformation complex
\[
(K(\stp_n\tmid_{\leq k}, \Poiss_n\{n\}\tmid_{\leq k}),d_0+d_1),
\]
see also \cite[Theorem 10.11]{FTW}.
We similarly define the truncated version $D_{n,g,k}$ of the genus $g$ piece so that
\[
(K(\stp_n\tmid_{\leq k}, \Poiss_n\{n\}\tmid_{\leq k}),d_0+d_1) \cong \prod_g D_{n,g,k},
\]

\subsection{Proof of Proposition \ref{prop:mainGC}}
We first note that in the case $m=n$ our formula \eqref{equ:d} for the homological degree simplifies to 
\begin{equation}\label{equ:dDgn}
d=
(n-2)g
-h-r
+3.
\end{equation}
We now try to proceed similarly to subsection \ref{sec:HGCproof}.
We consider the complex $D_{n,g}$ for fixed $g$.
It contains pieces of various $h,r$, with (repeating \eqref{equ:rbounds})
\begin{equation}\label{equ:Rbounds}
g+h\geq r\geq h\geq 1.
\end{equation}
The homological degree depends only on the number $h+r$, and we may consider $-h-r$ as the homological degree, up to an overall shift.
We want to study the homology map induced by the projection to the arity truncation
\begin{equation}\label{equ:DngDngk}
D_{n,g} \to D_{n,g,k}.
\end{equation}
Given the bounds on $h,r$, the above complexes can be depicted as follows, in $r,h$-space.
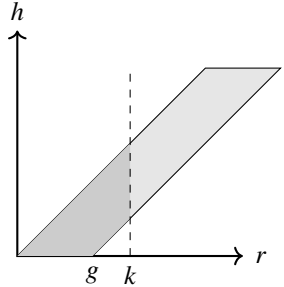
\[
\begin{tikzcd}
\draw[<->,thick] (0,3)--(0,0)--(3,0);
\node[above] at (0,3) {h};
\node[below] at (1,0) {g};
\node[right] at (3,0) {r};
\draw[fill=black!10] (0,0) -- (2.5,2.5) -- (3.5,2.5) -- (1,0) --(0,0);
\path[fill=black!20] (0,0) -- (1.5,1.5) -- (1.5,0.5) -- (1,0) --(0,0);
\draw[dashed] (1.5,0) node[below]{k} -- (1.5,2.5); 
\end{tikzcd}
\]
Here the shaded strip corresponds to $(r,h)$ such that pieces with that hair number and arity can be present in the complex $D_{n,g}$. Similarly the dark shaded piece corresponds to those pairs $(r,h)$ where the quotient $D_{n,g,k}$ can be nonzero. The vertical line at $r=k$ represents the truncation.
Note that in contrast to the hairy situation the truncation now is not just by the homological degree, but instead the truncation line cuts through various degree pieces.
This makes the effect of the truncation less transparent, and is the reason why the statement of Proposition \ref{prop:mainGC} is a bit weaker than that of Theorem \ref{thm:mainHGC}.

Concretely, consider the subspace $V_R$ of $D_{n,g}$ spanned by all pieces having $r+h=R$. This is the same as the subspace of homological degree $(n-2)g-R+3$. 
Then from $h\leq r$ we obtain $R=r+h\leq 2r$, or
\[
r\geq \frac R 2.
\]
Similarly, $h\geq 1$ gives $R\geq r+1$ and $h\geq r-g$ (see \eqref{equ:Rbounds}) yields $R\geq 2r-g$, or
\[
r\leq \frac{R+g}2.
\]
Hence the arities $r$ (potentially) present in $V_R$ are in the range
\[
r_{min}:= \left\lceil \frac R 2 \right\rceil \leq r \leq \min \left\{R-1,\left\lfloor \frac {R+g} 2 \right\rfloor  \right\} =: r_{max}
\]

In particular we see that the projection \eqref{equ:DngDngk} has the following properties:
\begin{itemize}
\item $V_R$ is sent to zero if $k< r_{min}$.
This is satisfied if
\[
R> 2k.
\]
\item $V_R$ is sent to itself if $k\geq r_{max}$. This is satisfied if
$k\geq R-1$ and $2k\geq R+g$, and hence if
\[
R\leq 2k - \min(k-1,g).
\]
\end{itemize}

Translating this back to homology, we see that the map
\[
H_d(D_{n,g}) \to H_d(D_{n,g,r})
\]
\begin{itemize}
	\item ... is zero if 
	$$d< (n-2)g +3-2k,$$
	 because the right-hand side $H_d(D_{n,g,r})$ is zero there.
	\item ... is an injection if 
	\[
	d\geq (n-2)g +3 
	-2k +\min(k-1,g)
	\]
	and an isomorphism (at least) if the inequality is strict.
\end{itemize}

In between these two degree bounds we may find (potentially) non-trivial pieces in the cokernel.
To arrive at the statements of Proposition \ref{prop:mainGC} we merely have to mind the additional degree shift in \eqref{equ:DngGC}

\hfill\qed

\subsection{Proof of Lemma \ref{lem:degboundGC}}
The Lemma has essentially been shown in \cite{Will}, in different degree conventions (cohomological instead of homological).
We quickly give a proof here for the sake of completeness.
The lower degree bound comes from combinatorial considerations on graphs. Concretely, in loop order $g\geq 2$ we may require all vertices to be at least trivalent, so that we have the inequality
\[
2e \geq 3v
\]
between the number of edges $e$ and the number of vertices $v$. Furthermore, we can express the loop number as $g=e-v+1$, and use this to eliminate $v=e-g+1$ in favor of $g$. In particular, our inequality becomes
\[
e \leq 3g-3.
\]
It follows that the degree (in $\GC_n^2$) satisfies the inequality
\[
d=(n-1)e-n(v-1)
=
(n-1)e-n(e-g)
=ng-e \geq (n-3)g+3.
\]

The upper degree bound for homology can be shown starting from the space $D_{n,g}$ above and \eqref{equ:DngGC}. To this end consider \eqref{equ:dDgn} and note that $r\geq 2$ and $h\geq 1$. For the degree in $D_{n,g}$ (still (ab)using the same letter $d$) we hence have
\[
d=(n-2)g-h-r+3\leq (n-2)g.
\]
Mind that by the degree shift in \eqref{equ:DngGC} we have to add $1$ to obtain a degree bound for $H(\GC_n^2)$.
The estimate can be slightly improved. Concretely, the $r=2$-part of the complex $D_{n,g}$ is easily written down explicitly. It is concentrated in loop orders $g\leq 2$. Hence for $g\geq 3$ we can use that $r\geq 3$ and the upper degree bound follows from \eqref{equ:dDgn}, where we have to mind the degree shift in \eqref{equ:DngGC}.


\end{document}